\documentclass[11pt, a4paper]{article}
\usepackage{amsmath}
\usepackage{amsfonts}
\usepackage{bm}
\newcommand\R{{\mathbb {R}}}
\newcommand\N{{\mathbb {N}}}
\newcommand{\Tm}{{\mathbb T}_m}

\newtheorem{Theorem}{Theorem}

\newtheorem{Lemma}[Theorem]{Lemma}

\newtheorem{Proposition}[Theorem]{Proposition}
\newtheorem{remark}[Theorem]{Remark}

\renewcommand{\P}{\ensuremath{\mathbb {P}}}
\newcommand{\E}{\mathbb{E}}

\newcommand{\Z}{ { \mathbb Z} }

\newcommand{\G}{{\mathcal G}}

\newcommand{\un}{\underline}
\newcommand\beq{\begin{equation}}
\newcommand\eeq{\end{equation}}

\def\cov{\mathop{\rm Cov}\limits}

\begin{document}

\title{An alternative to the coupling of Berkes-Liu-Wu
for strong approximations}

\author{Christophe Cuny$^a$, J\'er\^ome Dedecker$^b$ and Florence 
Merlev\`ede$^c$}

\maketitle

{\abstract{In this paper we propose an alternative to the coupling of Berkes, Liu and Wu \cite{BLW14}  to obtain  strong approximations for partial sums of dependent sequences.  The main tool is a new Rosenthal type inequality expressed in terms of the coupling coefficients. These coefficients are well suited  to some classes of   Markov chains or 
dynamical systems, but they also give new results for smooth functions of linear processes.}

\bigskip

\small{

$^a$ Universit\'e de la Nouvelle-Cal\'edonie, Institut de Sciences 
Exactes et Appliqu\'ees.

Email: christophe.cuny@univ-nc.nc \bigskip

$^b$ Universit\'e Paris Descartes, Sorbonne Paris Cit\'e,  Laboratoire MAP5 (UMR 8145).

Email: jerome.dedecker@parisdescartes.fr \bigskip

$^c$ Universit\'{e} Paris-Est, LAMA (UMR 8050), UPEM, CNRS, UPEC.

Email: florence.merlevede@u-pem.fr}

\section{Introduction}
Let $(\varepsilon_i)_{i \in {\mathbb Z}}$ be a sequence of
independent and identically distributed (iid) random variables,
and let $(X_n)_{n \in {\mathbb Z}}$ be a  strictly stationary sequence 
such that
\begin{equation}\label{firstseq}
X_n=f(\ldots,\varepsilon_0,  \ldots, \varepsilon_{n-1}, \varepsilon_n) \, ,
\end{equation}
for some real-valued measurable function $f$. 

In 2005, Wu \cite{Wu} introduced the so-called {\it physical 
dependence measure} defined in terms of  the following coupling:
let $\varepsilon_0'$ be distributed as $\varepsilon_0$ and 
independent of  $(\varepsilon_i)_{i \in {\mathbb Z}}$, and let
\beq \label{sequence0}
 \tilde X_n = 
 f(\ldots, \varepsilon_{-1},\varepsilon_0', \varepsilon_1, \ldots, \varepsilon_{n-1}, \varepsilon_n) \, .
\eeq
The physical  dependence coefficents in ${\mathbb L}^p$ (assuming that $\|X_0\|_p^p ={\mathbb E}(|X_0|^p) < \infty$) are 
then given by 
$$
\delta_p(n) = \|X_n - \tilde X_n\|_p \, .
$$
As pointed out by Wu, the coefficient $\delta_p(n) $ can be computed for a large variety of examples, including  iterated random functions and functions of linear processes. As we shall 
see, it is particularly easy to compute when $X_n$ is a smooth 
function of a linear process. 

Let $S_n=\sum_{k=1}^n X_k$. In a recent paper, Berkes, Liu and Wu \cite{BLW14}  use  the coupling defined above to prove the following 
strong approximation result: under an appropriate polynomial decay  of the coefficients  $\delta_p(n)$, the sequence $n^{-1} \E \big ( (S_n - n \E (X_1) )^2 \big )$ converges to $\sigma^2$ as $n \rightarrow \infty$ and, if $\sigma^2>0$, one can redefine $(X_n)_{n \geq 1}$ without changing its distribution on a (richer) probability space on which 
there exist iid random variables $(N_i)_{i \geq 1}$ with common distribution ${\mathcal N} (0, { \sigma}^2)$, such that,
\[
\left |\, S_n - n \E (X_1)- \sum_{i=1}^n N_i \, \right |=
 o\left (n^{1/p} \right ) \, \text{ ${\mathbb P}$-a.s.}
\]
The proof (from which we extract Proposition \ref{generalpropBLW},  Section \ref{Sec-IntroProof}) is based on a approximation by $m$-dependent sequences combined with an application of a deep result by Sakhanenko \cite{Sa06}.

This is a very important result, because it gives a full  extension 
of the Komlos, Major and Tusnady  (KMT) strong approximation \cite{KMT} for partial sums of iid random variables in ${\mathbb 
L}^p$   (for which $\delta_p(n)=0$ if $n\geq1$). 
 Most of the previous results in the dependent context were limitated to the rate $n^{1/4}$, because they were based on the Skorokhod representation theorem for martingales. As an exception, let us mention the paper 
\cite{MR}, where the rate $O(\log n)$ 
is reached for bounded observables of geometrically ergodic Markov chains.

In this paper, we follow the main steps of the proof of Berkes, Liu and Wu \cite{BLW14}, but we use a different coupling.  Let $(\varepsilon_i')_{i \in 
{\mathbb Z}}$ be an independent copy of  $(\varepsilon_i)_{i \in 
{\mathbb Z}}$
and let 
$$
X_n^*=f(\ldots, \varepsilon'_{-2}, \varepsilon'_{-1},\varepsilon_0', \varepsilon_1, \ldots, \varepsilon_{n-1}, \varepsilon_n) \, .
$$
Our coupling coefficient in ${\mathbb L}^p$ is then defined as
$$
\tilde \delta_p(n)= \|X_n-X_n^*\|_p \, 
$$
As $\delta_p$, this coefficient can be computed for a large class of 
examples (see Section \ref{examples}). 

Notice that $\E(X_n|\varepsilon_1,\ldots, \varepsilon_n)=\E(\tilde X_n|\varepsilon_1,\ldots, \varepsilon_n)$ a.s.,  and that 
$\|X_n-\E(X_n|\varepsilon_1,\ldots, \varepsilon_n)\|_p\le \tilde \delta_p(n)$, from which we easily deduce that
$\delta_p(n)\le 2 \tilde \delta_p(n)$.
Moreover, it  seems natural to think that, in many situations, the coefficient $\delta_p(n)$  should be much
smaller than $\delta'_p(n)$, because $X_n$ differs from $\tilde X_n$ by 
changing only the coordinate at point 0, while all the coordinates 
before time  0 have been changed in $X_n^*$. For Markov chains, 
however, the two coefficients should be of the same order (we
shall give in Section \ref{SecMR} an alternative definition of $\tilde \delta_p(n)$, which is more  adapted to the Markovian setting). 

Since $\delta_p(n) \leq 
 2 \tilde \delta_p(n)$, a reasonable question is then: what could be the  interest to deal with $\tilde \delta_p(n)$? The answer is simple :
the coupling $X_n^*$ is often easier to handle than 
$X_n'$ (because $X_n^*$ is by definition independent 
of the past  $\sigma$-algebra 
${\mathcal F}_0= \sigma(\varepsilon_i, i \leq 0)$) and 
we can develop specific tools involving the coefficent $\tilde \delta_p(n)$. 
In this paper, we shall prove and use a new Rosenthal-type
inequality (see Section \ref{SecRI}) expressed in terms of the coefficients $\tilde \delta_2$ and $\tilde \delta_p$. As a consequence, the conditions that  we impose on  $ \tilde \delta_p(n)$ are weaker 
than the corresponding conditions on $\delta_p(n)$ in the paper by 
Berkes, Liu and Wu. The two results are not comparable, 
but we shall obtain better conditions in all the cases 
where $\delta_p(n)$ and $\tilde \delta_p(n)$ are exactly
of the same order (for instance in the case of Markov chains). 

Let us present a simple example where our conditions are less restrictive than 
those of Berkes, Liu and Wu. Assume that
$$
X_n= g \left ( \sum_{i=0}^\infty a_i \varepsilon_{n-i}
\right )
$$
where $(a_i)_{i \geq 0} \in \ell_1$, and $(\varepsilon_i)_{i  \in 
{\mathbb Z}}$ is a sequence of iid random variables in 
${\mathbb L}^p$. Here $g$ is a continuous function
such that 
$$
|g(x)-g(y)| \leq c(|x-y|)
$$
where $c$ is a non-decreasing concave function and $c(0)=0$ ($c$ is then a concave majorant  of the modulus of 
continuity of $g$). In that case, using Lemma 5.1 in 
\cite{Dedecker},  it is easy to see that
$$
\delta_p(n) \leq \|c(|a_n(\varepsilon_0-\varepsilon'_0)|)\|_p \leq c(2\|\varepsilon_0\|_p |a_n|) \, ,
$$
and 
$$
\tilde \delta_p(n) \leq \left \|c\left ( \left |\sum_{i=n}^\infty a_i(\varepsilon_{n-i}-\varepsilon'_{n-i})\right | \right ) \right \|_p 
\leq c\left ( \left \|\sum_{i=n}^\infty a_i(\varepsilon_{n-i}-\varepsilon'_{n-i})\right \|_p \right )
\leq c\left ( C_p \|\varepsilon_0\|_p \sqrt{ \sum_{i=n}^\infty a_i^2} \right ) \, ,
$$
where we have used Burkholder's inequality for the 
last upper bound  (the positive constant $C_p$ depends only on $p$).
As expected, we see that the upper bound for $\delta_p(n)$ is smaller than the upper bound for $\delta'_p(n)$. 

Let us consider now the case where $|a_i|= O(i^{-a})$ 
for some $a>1$, and $c(x)\leq C |x|^\beta$ in a neighborhood of $0$ for some $\beta \in (0,1]$. 
In that case, the conditions  of Berkes, Liu and Wu on 
$\delta_p$
hold provided
$$
a > \frac 2 \beta \quad \text{for $p \in (2,4]$, and} 
\quad a> \frac{\tau(p)+1}{\beta} \quad 
\text{for $p>4$,}
$$
where 
$$
\tau(p):=\frac{(p-2)\sqrt{p^2+20p+4}\,+p^2-4}{8p}\, ,
$$
while our condition on $\tilde \delta_p$ are satisfied as soon
as
$$
a>\frac{\kappa(p)}{\beta} + \frac{1}{2} \, ,
$$
where, for $p>2$,
\begin{equation}\label{def-kappa}
\kappa(p):=\frac{(p-2)\sqrt{p^2+12p+4}\,+p^2+4p-4}{8p}\, .
\end{equation}
As one can see, our condition on $a$ is always less restrictive: for $p \in (2, 4]$ it suffices to notice that 
$\kappa$ is increasing and $\kappa(4)<1.4$. For $p>4$, 
it suffice to notice that $\kappa(p)<\tau(p) +0.5$. Note that, since 
$\kappa(p) \rightarrow 1/2$ as $p \rightarrow 2$, in the case where 
$\beta=1$ (Lipschitz observables), we only need $a>1$ and a moment
of order $b>2$ for $\varepsilon_0$ to get a strong approximation 
of order $o(n^{1/(2 + \epsilon)})$ for some $\epsilon>0$.

In addition to this example of  functions of linear processes, 
we shall apply our main results to some classes of Markov chains  or dynamical systems. The Markov chains we shall consider are not (or have no reasons to be) 
 irreducible, and some kind 
of regularity on the observables is required  (as in the previous example). We shall express these regularity conditions  in terms of the modulus of continuity (or ${\mathbb L}^p$-modulus of continuity) of the observables. These examples of Markov chains are different from the examples we considered in the previous paper \cite{CDM}, where we used the coefficient $\tilde \delta_1$. 
On the one hand, the coupling coefficient $\tilde \delta_1(n)$ can be computed  for  a larger
class of examples, but on the other hand  we need to impose a moment condition related to $p$ and to the decay rate of $\tilde \delta_1(n)$  to get the strong approximation 
with rate $o(n^{1/p})$.

\smallskip

In all the paper, we shall use the notation
$a_n \ll b_n$, which means that there exists a positive
constant $C$ not depending on $n$ such that
$a_n \leq C b_n$, for all positive integers
$n$.

\section{Main results}\label{SecMR}

Before giving our first main result, let us give the appropriate definition of the coefficient 
$\tilde \delta_p$ when $(X_n)_{n \in {\mathbb Z}}$ is a stationary sequence such that  
\beq\label{process}
X_n:=f(\varepsilon_n,\varepsilon_{n+1},\ldots)\, 
\eeq
for some measurable real-valued function $f$. 
This representation will play an important role in the application  to certain non-invertible dynamical systems.

Recall that  $(\varepsilon_i)_{i\ge 0}$ is a sequence of iid random variables,  and  that  $(\varepsilon_i')_{i\ge 0}$ is an 
independent copy of $(\varepsilon_i)_{i\ge 0}$. 
Define then $ \tilde X_{1,n}:=f(\varepsilon_1,\ldots , \varepsilon_{n}, 
\varepsilon'_{n+1}, \varepsilon'_{n+2}, \ldots)$. Then, for 
every $p\geq 1$, the  coefficient $\tilde \delta_p(n)$ is 
defined by:
\beq\label{coeff}
\tilde \delta_p (n):= \|X_{1}-\tilde X_{1,n}
\|_p\, .
\eeq
Recall also that, for any $p>2$,  the fonction $\kappa(p)$
has been defined in \eqref {def-kappa}.

\begin{Theorem} \label{KMTavecdelta'p_bis} Let $(X_n)_{n \in {\mathbb Z}}$ be a stationary sequence defined by either \eqref{firstseq} or  \eqref{process}, and assume that $X_0$ has a moment of order $p>2$. Assume in addition that there exists a positive constant $c$ such that for any $n \geq 1$, 
\beq \label{condondelta'p-bis}
{\tilde \delta}_{p} (n) 
\leq  c n^{-\gamma} \, ,
\eeq 
for some $\gamma> \kappa(p)$.
Let $S_n = \sum_{k=1}^{n} X_k$. Then $n^{-1} \E \big ( 
(S_n - n \E (Y_1) )^2 \big ) \rightarrow { \sigma}^2$ as $n \rightarrow \infty$ and one can redefine $(Y_n)_{n \geq 0}$ without changing its distribution on a (richer) probability space on which 
there exist iid random variables $(N_i)_{i \geq 1}$ with common distribution ${\mathcal N} (0, { \sigma}^2)$, such that,
\[
\left |\, S_n - n \E (X_1)- \sum_{i=1}^n N_i \, \right |= o\left (n^{1/p} \right ) \quad \text{ ${\mathbb P}$-a.s.}
\]
\end{Theorem}

\begin{remark}
Concerning the function $\kappa$, note that 
$\kappa(p) < (p+4)/4$ and that the function $p\rightarrow (p+4)/4$ is an asymptot of $\kappa$
as $p\rightarrow \infty$. 
\end{remark}

As quoted in the introduction, we shall now consider the case where the variables $X_n$ are  functions  of  random iterates. In that case the representation \eqref{firstseq} is not necessarily appropriate (nor even  easy to establish), and we need to define an appropriate coefficient $\delta'_p$ similar to $\tilde \delta_p$. 

Let $(\varepsilon_i)_{i \geq 1}$ be iid random variables with values in a measurable space $G$ and common distribution $\mu$. Let $W_0$ be a random variable with values in a measurable space $X$, independent of 
$(\varepsilon_i)_{i \geq 1}$ and let 
$F$ be a measurable function from $G\times X $ to $X$. For any $n \geq 1$, define
\beq \label{def-W}
W_{n} = F ( \varepsilon_{n}, W_{n-1} ) \, ,
\eeq
and assume that $(W_n)_{n \geq 1}$ has a stationary distribution $\nu$. 
Let now $h$ be a measurable function from $G\times X $ to ${\mathbb R}$ and define, for any $n \geq 1$,
\beq \label{functionofarandomiterates}
X_n = h ( \varepsilon_{n}, W_{n-1} ) \, .
\eeq
Then $(X_n)_{n \geq 1}$ is a stationary sequence with stationary distribution, say  $\pi$.  Let $( {\mathcal G}_{i} )_{i \in {\mathbb Z}}$ be the non-decreasing filtration defined as follows: for any $i < 0$, 
${\mathcal G}_{i} =\{ \emptyset, 
\Omega\}$, ${\mathcal G}_{0} = \sigma (W_0)$ and for any $i \geq 1$, $ 
{\mathcal G}_{i} = \sigma ( \varepsilon_{i }, \ldots,  \varepsilon_{ 1 }, W_0 ) $. It follows that for any $n \geq 1$, $X_n$ is ${\mathcal G}_{n}$-measurable. 

Let $W_0$ and $W_0^*$ be two random variables with 
law $\nu$, and such that $W_0^*$ is independent of
$(W_0, (\varepsilon_i)_{i \geq 1})$. For any $n\geq 1$, let 
$$
X_n^*=h(\varepsilon_n, W^*_{n-1}) \ \text{with} \ 
W_n^*=F(\varepsilon_n, W_{n-1}^*) \, .
$$
%Let also $(X_{n,x})_{n \geq 1}$ be the sequence of  %random variables 
%starting from $W_0 =x$.
We then define the coefficients  $(\delta_p'(n))_{n\ge 0}$ as follows
\begin{equation}\label{def-delta-prime}
\delta_p'(0) :=  \|X_1\|_p \,  \text{ and } \, \delta_p' (n) :
= \sup_{k \ge n} \| X_{k} - X_{k}^*\|_p  ,  \ n \geq 1 \, . 
\end{equation}
It is not difficult to see that, for any
positive integer $n$, 

\begin{equation}\label{simple-bound}
\delta_p' (n) \le 2\| X_{n} - X_{n}^*\|_p \, .
\end{equation}

For such functions  of random iterates, the following 
counterpart of Theorem \ref{KMTavecdelta'p_bis}
holds:

\begin{Theorem} \label{KMTavecdelta'p} Let $(X_n)_{n \geq 1}$ be a stationary sequence defined by \eqref{functionofarandomiterates} and assume that its stationary distribution $\pi$  has a moment of order $p>2$. Assume in addition that there exists a positive constant $c$ such that for any $n \geq 1$, 
\beq \label{condondelta'p}
{\delta}_{p}' (n) 
\leq   c n^{-\gamma} \, ,
\eeq 
for some $\gamma> \kappa(p)$.
Let $S_n = \sum_{k=1}^n X_k$. Then $n^{-1} \E \big ( (S_n - n \E (X_1) )^2 \big ) \rightarrow { \sigma}^2$ as $n \rightarrow \infty$ and one can redefine $(X_n)_{n \geq 1}$ without changing its distribution on a (richer) probability space on which 
there exist iid random variables $(N_i)_{i \geq 1}$ with common distribution ${\mathcal N} (0, { \sigma}^2)$, such that,
\[
\left |\, S_n - n \E (X_1)- \sum_{i=1}^n N_i \, \right|= o\left (n^{1/p} \right ) \, \text{ ${\mathbb P}$-a.s.}
\]
\end{Theorem}

\section{Applications}\label{examples}

\subsection{Applications to contracting iterated random functions.} 

We use the notations from the second part of Section \ref{SecMR}, with a Markov chain $(W_n)_{n \geq 1}$ defined by the recursive equation \eqref{def-W} and a sequence $(X_n)_{n \geq 1}$ defined by \eqref{functionofarandomiterates}. Assume that $X$ is equipped with a  metric $d$ and that it is endowed with the corresponding Borel $\sigma$-algebra.  Let us fix a ``base point" 
$x_0\in X$. For every $x\in X$, write $\chi(x):= 1+d(x_0,x)$.

Let us assume that there exists $C>0$, $\rho\in (0,1)$ and 
$\alpha\ge 1$, such that 

\beq
\label{start}\int_G \big(\chi(F(g,x_0)) \big)^\alpha \mu(dg) <\infty\, ,
\eeq
and 
\beq
\label{lip-2}\E\big(d(W_{n,x},W_{n,y})^\alpha\big)\le C \rho^n (
d(x,y))^\alpha\, ,
\eeq
\smallskip
where $W_{n,x}$ is the chain defined by \eqref{def-W} 
starting from $W_0=x$.

The next lemma is a combination of Theorem 2 and 
Lemma 1 of Shao and Wu \cite{SW}.

\begin{Lemma}\label{lemSW}
Assume that \eqref{start} and \eqref{lip-2} hold for some $\alpha\ge 1$, 
$C>0$ and $\rho\in (0,1)$. Then the Markov chain $(W_n)_{n\in \N}$ admits 
a stationary distribution $\nu$ such that 
$\int_X (\chi(x))^\alpha \nu(dx) <\infty$. Moreover,  there exist $C(
\alpha)>0$ and $\rho(\alpha)\in (0,1)$, such that for every $n\ge 1$, 
$$
\iint \E\big( (d(W_{n,x},W_{n,y}))^\alpha\big) \nu(dx) \nu(dy) \le 
C(\alpha) (\rho(\alpha))^n\, .
$$
\end{Lemma}

From now, the sequence $(X_n)_{n\geq 1}$ is defined 
by  \eqref{functionofarandomiterates}, where
the chain $(W_n)_{n\geq 0}$ is strictly stationary, 
with stationary distribution $\nu$.

We shall say that a function $h \, :\, G\times X\to \R$ satisfies the assumption $H_{s,t}$, for some $s,t \geq 0$ if there exist non-negative 
 functions $\eta$ and $\tilde \eta$ and a non-decreasing function $\beta\, :\, 
[0,1]\to [0,+\infty)$ such that, for every $(g,x)\in 
G \times X$, 

\beq \label{moment-2}
\big|h(g,x)\big|\le \eta(g) (\chi(x))^s\, ,
\eeq
and for  every $u\in (0,1]$,

\beq\label{lip-3}
\sup_{y\in X\, :\, d(x,y)\le u} \big|h(g,x)-h(g,y)\big|
\le \beta(u) \tilde \eta(g) (\chi(x))^{t}\, .
\eeq

\begin{Lemma}\label{contract}
Let $p\ge 1$.
Assume that \eqref{start} and \eqref{lip-2} hold for some $\alpha\ge 1$. 
Assume that \eqref{moment-2} and \eqref{lip-3} hold for 
some $0\le s<  \alpha/p$, $0\le t\le \alpha/p$ and some $\eta,\tilde \eta$ such that 
$\int_G (\eta(g))^p\mu(dg) <\infty$ and $\int_G(\tilde \eta(g))^p \mu(dg)
<\infty$. Then, there exist $0<\omega_1,\omega_2<1$ and  $C>0$, such that 
\beq
\delta_p'(n)\le C(\beta(\omega_1^n) + \omega_2^{n})\qquad \forall n\ge 1\, .
\eeq
\end{Lemma}

\noindent {\bf Proof.} 
Let $\varepsilon >0$, and let $u_n(x,y):= \big|h(\varepsilon_{n+1}, 
W_{n,x})-h(\varepsilon_{n+1}, W_{n,y})\big|$. 
Writing, 
\begin{equation}
u_n^p(x,y)=u_n^p(x,y) {\bf 1}_{\{d(W_{n,x},W_{n,y})\le \rho^{\varepsilon n}\}} 
+ u_n^p(x,y){\bf 1}_{\{ d(W_{n,x},W_{n,y})>\rho^{\varepsilon n} \}} \, ,
\end{equation}
we obtain the upper bound 
$$
\iint \E(u_n^p(x,y))\,\nu(dx)\nu(dy) \le I_n+II_n\, .
$$

Clearly,

$$
I_n \le \left (\beta(\rho^{\varepsilon n}) \right )^p \left(\int_G (\tilde \eta (g) )^p\mu(d g)
\right) \left ( \int_X( \chi(x))^{pt} \nu(dx)\right)\, .
$$

Moreover, using H\"older's inequality and  %$\int_X(\chi(x))^\alpha \nu(dx) <\infty$ by 
Lemma 
\ref{lemSW},  we have

\begin{multline*}
II_n \\ 
\le \frac{2^{p-1}}{\rho^{\varepsilon (\alpha-sp)n}} \left (\int_G \eta^p(g)\mu(dg)
\right )\iint \E\left[ \left(\chi^{sp}(W_{n,x})+\chi^{sp}(W_{n,y})\right)
(d(W_{n,x},W_{n,y}))^{\alpha-sp}\right]\, \nu(dx)\nu(dy)\\
\le \frac{C}{\rho^{\varepsilon (\alpha-sp)n}} \left (\iint \E(d(W_{n,x},W_{n,y}))^{\alpha}\, \nu(dx)\nu(dy) 
\right) ^{1-sp/\alpha}\le \tilde C 
\left (\frac{(\rho(\alpha))^{1/\alpha}}{\rho^{\varepsilon}}\right)^{(\alpha-sp)n}\,  .
\end{multline*}
The desired bound follows by taking $\varepsilon$ small enough. 
\hfill $\square$

\medskip

\begin{Proposition} Let $1<p<\alpha$. Assume that there exist $C>0$ and 
$\rho\in (0,1)$ such that \eqref{start} and \eqref{lip-2} hold. 
Let $0\le s<\alpha/p$ and $0\le t\le p/\alpha$. Let $h$ 
satisty $H_{s,t}$. Assume that $\int_G (\tilde \eta(g))^p\mu(dg)<\infty$, 
$\int_G (\eta(g))^{p} \mu (dg) < \infty$ and that 
$\beta(2^{-n})=O(n^{-\gamma})$, with $\gamma> \kappa(p)$. Then, the conclusion of Theorem \ref{KMTavecdelta'p} holds.
\end{Proposition}
\noindent {\bf Proof.} 
Starting from Theorem \ref{KMTavecdelta'p} and Lemma 
\ref{contract},  it is enough to prove that our assumption implies that, for any $a\in 
(0,1)$, $\beta(a^{n})=O(n^{-\gamma})$. 
Too see this, we note that there exists an integer $\ell \ge 1$ such that 
$a ^\ell \le 1/2$. Let $n\ge 2\ell $. 
Notice that $n/(2\ell)\le [n/\ell]\le n/\ell$. Since $\beta$ is non-decreasing, we have 
$$
\beta(a^n)\le \beta(2^{-[n/\ell]})\le C([n/\ell])^{-\gamma}\le C(2\ell /n)^\gamma\, ,
$$
and the result follows. \hfill $\square$

\subsection{Applications to dilating endomorphisms of the torus}
\label{torus}

Let $A$ be an $m\times m$ matrix with integral entries. Then, 
$A$ induces a  transformation $\theta_A$ of the $m$-dimensional torus ${\mathbb T}_m:=\R^m/\Z^m$ preserving the Haar measure 
$\lambda$.

Assume that $A$ is dilating, i.e. that all its eigenvalues have modulus 
strictly greater than one. Let $\Gamma$ be a system of representative of 
$\Z^m/A\Z^m$.
Then, $\theta_A$ admits a Perron-Frobenius operator 
$P_A$ given by 

\begin{equation}\label{perron-frobenius}
P_Af(x)=\frac1{N}\sum_{\gamma\in \Gamma} 
f(A^{-1}x+A^{-1}\gamma)\, ,
\end{equation}
for every continuous function $f$ on $\Tm$, where $N=|{\rm det }\, A|=\#\Gamma$.

Since $P_A$ is markovian, there exists a Markov chain with state space 
${\mathbb T}_m$ admitting $\lambda$ as stationary distribution. 
This Markov chain may be realized as follows:
let $W_0$ be a random variable taking values in $\Tm$ and 
$(\varepsilon_i)_{i\ge 1}$ be iid variables uniformly distributed on $\Gamma$ 
and independent of $W_0$. For every $n\ge 1$, define 
$W_n:= A^{-1}W_{n-1}+A^{-1}\varepsilon_n$. 
Denote by $(W_{n,x})_{n\ge 0}$ 
 the Markov chain starting at $x\in \Tm$.
 
 Let $h$ be some measurable function from $\Tm$ to ${\mathbb R}$,  and 
 let $X_n=h(W_n)$ where $W_0$ has distribution $\lambda$. 
 Let also $X_{n,x}=h(W_{n,x})$.

For every $p\ge 1$ and every $f\in {\mathbb L}^p(\lambda)$  the  ${\mathbb L}^p$-modulus of continuity of $f$ 
is given by
$$
\omega_{p,f}(\delta):= \sup_{|x|\le \delta}\|f(\cdot + x)-f\|_p\quad \forall \ 0\le \delta\le 1\, ,
$$
where $|\cdot|$ stands for the euclidean norm.

\begin{Lemma}\label{lem-matrix}
Let $p\geq 1$ and  $h\in {\mathbb L}^p(\lambda)$. 
%For every 
%$n\in \N$ and every $x\in \Tm$, set $X_{n,x}:=h(W_{n,x})$.  
The following upper bound holds:
$$
\left( \iint \, \E(|X_{n,x}-X_{n,y}|^p)\, \lambda(dx)\lambda(y)\, \right )^{1/p}\le 2^{m/p} \omega_{p,h}\left(\Delta\left (A^{-n}([0,1]^m)\right ) \right) \, ,
$$ 
where $\Delta\left (A^{-n}([0,1]^m)\right )$ stands for the diameter of 
$A^{-n}([0,1]^m)$. Consequently (using \eqref{simple-bound}),  the coefficients $\delta'_p(n)$ of the stationary sequence $(X_n)_{n \in {\mathbb  Z}}$ satisfy 
$$
\delta'_p(n) \leq 2^{m/p+1} \omega_{p,h}\left(\Delta\left (A^{-n}([0,1]^m)\right ) \right) \, .
$$
\end{Lemma}
\noindent {\bf Proof.} 
We start by some preliminary considerations. Iterating the recursive equation
$W_{n,x}= A^{-1}W_{n-1,x}+A^{-1}\varepsilon_n$, we get that
$$
W_{n,x}= A^{-n} x + \sum_{i=1}^n A^{-i} \varepsilon_{n-i +1} \, .
$$
Note that the random variable $W_{n,x}$ has the same distribution as 
$Y_{1,x}$, where $Y_{1,x}$ is 
the first iteration of the Markov chain starting at $x$ with 
transition $P^n_A=P_{A^n}$. As explained at the beginning of this section, this 
may be realized as 
$$
Y_{1,x}= A^{-n} x + A^{-n} \xi_1 \, ,
$$
where $\xi_1$  is uniformly distributed 
over $\Gamma_n$ (a system of representative of  $\Z^m/A^n\Z^m$).
Let $Z_{1,x}=h(Y_{1,x})$ It follows that 
$$
\iint \, \E(|X_{n,x}-X_{n,y}|^p)\, \lambda(dx)\lambda(y)=
\iint \, \E(|Z_{1,x}-Z_{1,y}|^p)\, \lambda(dx)\lambda(y) \, .
$$

From this last equality, we see  that it suffices to prove Lemma 
\ref{lem-matrix} for 
$n=1$, the general case then follows by considering $A^n$ rather than 
$A$. 

We refer to \cite{CHR} for the results that we need about 
tiling. There exists a unique compact set $K\subset \R^m$, 
such that 
\begin{equation}\label{self-affine}
K=\cup_{\gamma\in \Gamma}(A^{-1}K+A^{-1}\gamma)
\end{equation}
 and an integer 
$q\ge 1$ such that 
$$
\sum_{\un n\in \Z^m} {\bf 1}_{K+\un n} =q \qquad \mbox{$\lambda$-almost everywhere.}
$$
%with $\lambda$ the Lebesgue measure on $\R^d$.
Moreover, for every $\gamma,\gamma'\in \Gamma$ with $\gamma\neq \gamma'$, 
$\lambda\big( (A^{-1}K+A^{-1}\gamma)\cap (A^{-1}K+A^{-1}\gamma')\big)=0$. 
Using that
$$
{\bf 1}_K= \sum_{\un n\in \Z^m} {\bf 1}_{(K\cap ([0,1]^m-\un n)} = \sum_{\un n\in \Z^m} {\bf 1}_{((K+\un n)\cap[0,1]^m)-\un n}\qquad \mbox{$\lambda$-almost everywhere}\, ,
$$
we then infer that for every 
$\Z^m$-periodic locally integrable function $g$ on $
\R^m$, 
\begin{equation}\label{per}
\int_K g\, d\lambda= \int_{\R^m} \Big (\sum_{\un n \in \Z^m}{\bf 1}_{((K+\un n)\cap[0,1]^m)-\un n}\Big )g\, d\lambda
 = \sum_{\un n \in \Z^m} 
 \int_{(K+\un n)\cap[0,1]^m}g\, d\lambda= q \int_{\Tm}g \,d\lambda\, .
\end{equation}

Let $h\in {\mathbb L}^p(\lambda)$ (we identify $h$ with a $\Z^m$-periodic function 
on $\R^m$). We have 
\begin{multline*}
\iint \, \E(|X_{1,x}-X_{1,y}|^p)\, \lambda(dx)\lambda(y)
=\iint \, \E(|h(A^{-1}x+A^{-1}\varepsilon_1)-h(A^{-1}y+A^{-1}\varepsilon_1)|^p)\, \lambda(dx)\lambda(dy)\\
=\int_{\Tm}   \left (\int_{\Tm-x} \E(|h(A^{-1}x+A^{-1}\varepsilon_1)-h(A^{-1}(x+y)+A^{-1}\varepsilon_1)|^p)\, \lambda(dy) \right )
\lambda(dx)\\
=\int_{[-1,1]^m}   \left (\int_{(\Tm-y)\cap\Tm} \E(|h(A^{-1}x+A^{-1}\varepsilon_1)-h(A^{-1}(x+y)+A^{-1}\varepsilon_1)|^p)\, \lambda(dx) \right ) 
\lambda(dy) \,  .
\end{multline*}
Set 
\begin{multline*}
\psi_y(x):=\E(|h(A^{-1}x+A^{-1}\varepsilon_1)-h(A^{-1}(x+y)+A^{-1}\varepsilon_1)|^p)\\ =\frac1N \sum_{\gamma\in \Gamma }|h(A^{-1}x+A^{-1}\gamma)-h(A^{-1}(x+y)+A^{-1}\gamma)|^p\, .
\end{multline*}
Notice that $\psi_y$ is $\Z^m$-periodic. Hence, using \eqref{per} 
and \eqref{self-affine}, we have 
\begin{multline*}
\int_{(\Tm-y)\cap\Tm} \psi_y(x)\lambda(dx)\le \int_{\Tm} \psi_y(x)\lambda(dx)
= \frac1q \int_K \psi_y(x)\lambda(dx) \\
=\frac1q \int_{\cup_{\gamma\in \Gamma}(A^{-1}K+A^{-1}\gamma)}|h(x)-h(x+A^{-1}y)|^p\lambda(dx)= \frac1q \int_{K}|h(x)-h(x+A^{-1}y)|^p\lambda(dx)\\
=\int_{\Tm}|h(x)-h(x+A^{-1}y)|^p\lambda(dx)\, ,
\end{multline*}
and the result follows. \hfill $\square$

\medskip

We shall now explain how to obtain  the strong approximation result with rate $o(n^{1/p})$ for the partial sums of the process 
 $(h\circ \theta_A^n)_{n\in \N}$ for $h\in {\mathbb L}^p(\lambda)$. 
Let $(\varepsilon_n)_{n\ge 0}$ be a sequence of iid variables uniformly distributed on $\Gamma$. We define a probability $\nu$ on ${\mathbb T}_m$ 
by setting, for every $f\in C([0,1])$,

$$
\int_{{\mathbb T}_m} f\, d\nu := \E\left (f \left(\sum_{k\ge 0}A^{-k-1}\varepsilon_k\right)\, \right)\, .
$$
By construction, $\nu$ is $P_A$-invariant. Since $A$ is dilating, the only $P_A$-invariant probability on ${\mathbb T}_m$ is $\lambda$.

Define $Z_0:= \sum_{k\ge 0}A^{-k-1}\varepsilon_k$ and for 
every $n\ge 1$, (with equality in $\Tm$)
$$Z_n:= A^n Z_0 =\sum_{k\ge 0}A^{n-k-1}\varepsilon_k
=\sum_{k\ge n}A^{n-k-1}\varepsilon_k=\sum_{k\ge 0}A^{-k-1}\varepsilon_{k+n}
\, .
$$
Notice that for any $h\in {\mathbb L}^p(\lambda)$ the processes 
$(h\circ \theta_A^n)_{n\ge 0}$ (under $\lambda$) and 
$(Y_n)_{n\in \N}:=(h(Z_n))_{n\ge 0}$ (under $\P$) have the same distribution. 

Let $\tilde \delta_p(n)$ be the coefficients associated with $(Y_n)_{n\ge 0}$ as in \eqref{coeff}. The computations done in the proof of Lemma 
\ref{lem-matrix} yield to the following bound 

\begin{equation}\label{bound-matrix}
\tilde \delta_p(n) \leq 2^{m/p+1} \omega_{p,h}\left(\Delta\left (A^{-n}([0,1]^m)\right ) \right)\, .
\end{equation}

%Hence we can also apply our strong invariance principles. 

\medskip

As a consequence of Lemma \ref{lem-matrix} and of \eqref{bound-matrix}, Theorem 
\ref{KMTavecdelta'p_bis} (applied to 
$(h(Z_n))_{n\ge 0}$) or Theorem \ref{KMTavecdelta'p} (applied to 
$(h(W_n))_{n\ge 0}$), lead to the following proposition:

\begin{Proposition} \label{KMTformatrix-piecewise} Let $p>2$
and let $\kappa(p)$ be defined in \eqref{def-kappa}. Let $h\in {\mathbb L}^p(\lambda)$ be such that $\omega_{p,h}(2^{-n})\le O(n^{-\gamma})$
for some $\gamma> \kappa(p)$.  Assume that, with the above notations, 
$S_n=X_1+\cdots +X_n=h(W_1)+\cdots + h(W_n)$, or 
$S_n =h(Z_1)+\cdots +h(Z_n)$. Then $n^{-1} \E \big ( (S_n - n \lambda(h) )^2\big ) \rightarrow \sigma^2$ as $n \rightarrow \infty$ and for every (fixed) $x\in[0,1]$, one can redefine $(S_{n})_{n \geq 1}$ without changing its distribution on a (richer) probability space on which 
there exist iid random variables $(N_i)_{i \geq 1}$ with common distribution ${\mathcal N} (0, \sigma^2)$, such that,
\[
\left | \, S_{n} - n \lambda(h) - \sum_{i=1}^n N_i 
\, \right | = o\left (n^{1/p} \right ) \, 
 \text{ ${\mathbb P}$-a.s.}
\]
\end{Proposition}

\begin{remark}
Alternatively, one can also apply Theorem 2 in \cite{CDM}, by using the upper bound on $\delta'_1(n)$  given in  Lemma \ref{lem-matrix}. For instance, if $h$ is bounded
and such that $\sum_{n>0} n^{p-2} \omega_{1,h}(2^{-n}) < \infty$, then the conclusion of Proposition 
\ref{KMTformatrix-piecewise}
holds. If $m=1$ (for instance for the transformation 
$\theta(x)=2x-[2x]$), this implies that, for $BV$-observables,  the strong approximation   holds with the 
rate $o(n^{1/p})$ for any $p>2$. 
\end{remark}

\subsection{Applications to dilating piecewise 
affine maps}

Let $K$ be a countable set, with $|K|\ge 2$, and $(I_k)_{k\in K}$ be a collection of 
disjoint open subintervals of $[0,1]$ such that $\bigcup_{k\in K}
\overline I_k
=[0,1]$ (it is possible to have several accumulation points).
Notice that  $\sum_{k\in K} \lambda(I_k)=1$, where $\lambda$ stands for the Lebesgue measure on $[0,1]$.

Let $T\, :\, [0,1]\to [0,1]$ be a map such that $T_{|I_k}$ is affine 
and onto $(0,1)$, so that $T_{|I_k}$ extends in a trivial way to 
an affine  map from $\overline I_k$ onto $[0,1]$, that we still denote by 
$T_{|I_k}$. The values of $T$ on $[0,1]\backslash \bigcup_{k\in K}I_k$ will be irrelevant in the sequel. 
For every $k\in K$, denote by $s_k$ the inverse of $T_{|I_k}$ from 
$(0,1)$ onto $I_k$. There exist reals $\alpha_k$ and $\beta_k$, such that 
for every $x\in (0,1)$, $s_k(x)=\alpha_k x +\beta_k$ (hence 
$T_{|I_k}(u)= (u-\beta_k)/\alpha_k$). Then, $|\alpha_k|=\lambda(I_k)$.

Such a map $T$ admits a Perron-Frobenius operator $P$ defined by 
$$
Pf(x):= \sum_{k\in K}|\alpha_k| f(\alpha_kx+\beta_k) \, ,
$$
for every continuous function $f$ on $[0,1]$.

Since $P$ is Markovian and leaves  $\lambda$ invariant, 
there exists a Markov chain with state space $[0,1]$ admitting $\lambda$ as stationary distribution. Since $|K|\ge 2$ then $0<|\alpha_k|<1$ for every 
$k\in K$ and one may easily prove  that $\lambda$ is the only $P$-invariant measure on $[0,1]$.

The above Markov chain may be realized as follows. Let $W_0$ be a random variable taking values in $[0,1]$. Let $(\varepsilon_i)_{i\ge 1}$ be iid 
random variables independent of $W_0$, taking values in $K$, 
such that $\P(\varepsilon_1=k)=|\alpha_k|$ for every $k\in K$. For every 
$n\ge 1$, set $W_n := s_{\varepsilon_n}(W_{n-1})$ and denote by 
$(W_{n,x})_{n\ge 0}$ the Markov chain starting from $x\in [0,1]$. 
Notice that for every $n\ge 1$ and every $x\in [0,1]$, 
$$
W_{n,x}=s_{\varepsilon_n}\circ \cdots \circ s_{\varepsilon_1}(0)
+\alpha_{\varepsilon_n}\ldots \alpha_{\varepsilon_1} x:= B_n +A_n x\, .
$$

Let $h$ be some measurable function from $[0,1]$ to ${\mathbb R}$,  and 
 let $X_n=h(W_n)$ where $W_0$ has distribution $\lambda$. 
 Let also $X_{n,x}=h(W_{n,x})$. 

For every $f\in C([0,1])$ define 
$$
\omega_{\infty,f}(\delta)=\sup_{x,y\in [0,1], \,|x-y|\le \delta}
|f(x)-f(y)|\, ,\quad \forall \delta\in [0,1]\, .
$$
Define also 
$$
\delta_\infty(n):= \sup_{x,y\in [0,1]}\E|X_{n,x}-X_{n,y}|\, .
$$

%where $|\cdot|$ stands for the euclidean norm.

\begin{Lemma}\label{lem-affine}
Let $h\in C([0,1])$, and let 
%For every 
%$n\in \N$ and every $x\in [0,1]$, set $X_{n,x}:=h(W_{n,x})$. 
$\bar \alpha:=\max_{k\in K}|\alpha_k|$. For 
 every integer $n\ge 1$, we have
$$
\sup_{x,y\in [0,1]} |X_{n,x}-X_{n,y}|\le 2\omega_{\infty, h}(\bar \alpha^n
) \, .
$$ 
In particular for every $n\ge 1$, 
$
    \delta_p'(n)\le 2 \omega_{\infty, h}(\bar \alpha^n
)$ for any $p\ge 1$, and 
$   \delta_\infty(n) \le 2 \omega_{\infty, h}(\bar \alpha^n)$.
\end{Lemma}
\noindent {\bf Proof.}  For every $x\in [0,1]$, we have 
\begin{multline*}
 |X_{n,x}-X_{n,y}|
=  |h(B_n+A_n x)-h(Z_n+A_n y)| \\ \le |h(B_n+A_n x)-h(B_n)|+|h(B_n+A_n y)-h(B_n)|\, ,
\end{multline*}
and the result follows. \hfill $\square$

\begin{remark} When $K=\{1,\ldots, r\}$ and 
$|\alpha_1|=\cdots =|\alpha_r|$, similar computations as those done in  Section \ref{torus} allow to control $(\delta_p'(n))_{n\in \N}$ 
thanks to the ${\mathbb L}^p$-modulus of continuity. Notice that, actually, 
the case where $\alpha_1=\cdots = \alpha_r$ is included in section 
\ref{torus} (taking $m=1$).
\end{remark}

\smallskip

As in the previous subsection, let us  also consider the process 
$(h\circ T^n)_{n\in \N}$. Let $(\varepsilon_n)_{n\in \N}$ 
be iid random variables taking values in $K$ such that 
$\P(\varepsilon_1=k)=|\alpha_k|$ for every $k\in K$. 
For every $n\in \N$, set $Z_n:= \sum_{\ell\in \N} s_{\varepsilon_n} 
\circ\cdots \circ s_{\varepsilon_{\ell+n}}(0)=\sum_{\ell\in \N} 
\big(\prod_{j=0}^{\ell-1}\alpha _{\varepsilon_{j+n}}\big)\,
 \beta_{\varepsilon_{\ell+n}}$. Then, $(Z_n)_{n\in\N}$ is identically 
 distributed and the common law is invariant by $P$, so it is the Lebesgue measure on $[0,1]$. Moreover, one can see that $Z_n=T^n Z_0$ for every $n\in \N$. 
 Hence, for every $h\in C([0,1])$, the processes $(h\circ T^n)_{n\in \N}$ 
 (under $\lambda$) and $(Y_n)_{n\in \N}:=(h(Z_n))_{n\in \N}$ (under $\P$) 
 have the same distribution. As above the following upper bound clearly holds
 
 \begin{equation}\label{bound-affine}
 \tilde \delta_\infty(n)\le  2 \omega_{\infty, h}(\bar \alpha^n)\, .
 \end{equation}

\begin{Proposition} \label{KMTforMarkov-piecewise} Let $p>2$,
$\kappa(p)$ be defined by \eqref{def-kappa},  and  $h\in C([0,1])$. Let
$S_n=X_1+\cdots + X_n=h(W_1)+ \cdots +h(W_n)$, or $S_n=h(Z_1)+\cdots +h(Z_n)$.
Assume that   $\omega_{\infty,h}(2^{-n})\le O(n^{-\gamma})$
for some $\gamma $  such that
$$
\gamma > \frac{p-1}{2} \ \text{if $p \in (2,3]$, and  
$\gamma > \kappa(p)$  if $p>3$.}
$$ Then $n^{-1} \E \big ( (S_n - n \lambda(h) )^2\big ) \rightarrow \sigma^2$ as $n \rightarrow \infty$ and for every (fixed) $x\in[0,1]$, one can redefine $(S_{n})_{n \geq 1}$ without changing its distribution on a (richer) probability space on which 
there exist iid random variables $(N_i)_{i \geq 1}$ with common distribution ${\mathcal N} (0, \sigma^2)$, such that,
\[
 \left |\,  S_{n} - n \lambda(h) - \sum_{i=1}^n N_i  \,
 \right  |= o\left (n^{1/p} \right ) \, 
 \text{ ${\mathbb P}$-a.s.}
\]
\end{Proposition}

\noindent {\bf Proof.} It suffices to combine Lemma \ref{lem-affine} 
(or \eqref{bound-affine})
with either our Theorem \ref{KMTavecdelta'p} (or our Theorem 
\ref{KMTavecdelta'p_bis}) for $p>3$ or  Theorem 1 of \cite{CDM} for $p \in (2,3]$ (which applies also to processes as defined by \eqref{process}).
\hfill $\square$

\section{Proof of the results} 
\label{Sec-IntroProof}

As in \cite{CDM}, the proof is based on  a general proposition 
that can be established by combining the arguments given in the paper by Berkes, Liu and Wu \cite{BLW14}. Let us now recall this   proposition: it applies to a strictly stationary sequence $(X_k)_{k \geq 1}$  of real-valued random variables in 
${\mathbb L}^p$ ($p>2$) that can be well approximated by a sequence of $m$-dependent random random variables, with the help  of an auxiliary sequence of iid random variables $( \varepsilon_i)_{i \geq 0}$.  Let $(M_k)_{k \geq 1}$ be a sequence of positive real 
numbers and define 
\beq \label{defphik}
\varphi_k (x) = (x \wedge M_k) \vee (-M_k)  \text{ and } g_k (x) = x-\varphi_k (x)  \, .
\eeq
Then, define
\beq \label{defXkj}
{ X}_{k,j} =\varphi_k (X_j)  - \E \varphi_k (X_j) \, \mbox{ and } \, 
W_{k, \ell} = \sum_{i=1+3^{k-1}}^{\ell +  3^{k-1}} { X}_{k,i}  \, .
\eeq
Let  now $(m_k)_{k \geq 1}$ be a non-decreasing sequence of positive integers such that $m_k =o ( 3^k) $, as $k \rightarrow \infty$,  and define 
\beq \label{defXkjtilde}
{\tilde X}_{k,j} = \E \big ( \varphi_k (X_j)  | \varepsilon_j, \varepsilon_{j-1}, \ldots, \varepsilon_{j-m_k}\big ) - \E \varphi_k (X_j) \text{ for any  $j \geq m_k + 1$} \, 
\text{ and } \, {\widetilde W}_{k, \ell} = \sum_{i=1+3^{k-1}}^{\ell +  3^{k-1}} {\tilde X}_{k,i} \, .
\eeq
Finally, set $k_0 := \inf \{ k \geq 1 \, : \,  m_k \leq 2^{-1} 3^{k-2} \}$ and define
\beq \label{defnuk}
\nu_k = m_k^{-1}  \big \{ \E ({\widetilde W}^2_{k, m_k} )  + 2 \E ({\widetilde W}_{k, m_k} ({\widetilde W}_{k, 2m_k}  -{\widetilde W}_{k, m_k} ))\big \} \, .
%(:= \Vert A_{k,1} \Vert_2^2 /(3m_k)) 
\eeq

\begin{Proposition}[Berkes, Liu and Wu \cite{BLW14}] \label{generalpropBLW} 
Let $p>2$. Assume that we can find a sequence of positive reals $(M_k)_{k \geq 1}$,  a  non-decreasing sequence of positive integers $(m_k)_{k \geq 1}$   such that $m_k =o ( 3^{2k/p} k^{-1}) $ as $k \rightarrow \infty$, in such a way that the following conditions are satisfied:
\beq \label{condtruncature}
\sum_{k \geq 1} 3^{k (p-1)/p}  \E (|g_k (X_1)|)< \infty \, , 
% \sum_{k \geq 1} 3^{k (p-1)/p}  \E (|X_1|{\mathbf 1}_{|X_1| >M_k})< \infty
\eeq
there exists  $\alpha \geq 1$ such that 
\beq \label{condmkdependence}
\sum_{k \geq k_0} 3^{- \alpha k/p} \left \Vert\max_{1 \leq \ell \leq 3^k - 3^{k-1}} \left |   { W}_{k, \ell} -  {\widetilde W}_{k, \ell} \right | \right \Vert_{\alpha}^{\alpha}  < \infty \, , 
\eeq
and there exists $ r \in ]2, \infty [$ such that 
\beq \label{condtoapplySakhanenko}
\sum_{k \geq k_0} \frac { 3^{k} }{  3^{kr /p}m_k}   \E \left (  \max_{1 \leq \ell \leq  3m_k }  \left |   {\widetilde W}_{k, \ell}    \right |^r \right ) < \infty \, .
\eeq
Assume in addition that 
\beq \label{series-variance} 
\mbox{the series }\quad  \sigma^2 =  {\rm Var} (X_1^2) + 2 \sum_{i \geq 1} {\rm Cov} (X_1, X_{i+1} )\quad 
\mbox{ converge,}
\eeq
% and that 
%\beq \label{defsigma2asalimit}
%\lim_{n \rightarrow \infty} \frac{{\rm Var} (S_n)}{n} = \sigma^2 \, , 
%\eeq
and
\beq \label{cond2v_k}
3^k (  \nu_k^{1/2}  - \sigma)^2 = o(3^{2k/p} (\log k)^{-1})  \, , \, \mbox{ as $k \rightarrow \infty$}\, .  
\eeq
Then, one can redefine $(X_n)_{n \geq 1}$ without changing its distribution on a (richer) probability space on which 
there exist iid random variables $(N_i)_{i \geq 1}$ with common distribution ${\mathcal N} (0, \sigma^2)$, such that,
\beq \label{conclusiongeneralprop}
 \left | \, S_n - n \E (X_1)- \sum_{i=1}^n N_i \, \right | = o\left (n^{1/p} \right ) \, \text{ ${\mathbb P}$-a.s.}
\eeq
\end{Proposition}

\smallskip

Theorem \ref{KMTavecdelta'p} is a 
consequence of  the next proposition, 
whose proof follows from Proposition \ref{generalpropBLW}.
This proposition applies to the stationary 
sequence $(X_n)_{n \geq 1}$ defined 
by \eqref{functionofarandomiterates} and 
the conditions are expressed in terms of the coefficient $\delta'_p$ defined in the second part of Section \ref{SecMR}. 
In what follows, all the proofs will be
 written with the help of that coefficient, but the arguments are exactly 
the same for the sequence defined by 
\eqref{firstseq} or \eqref{process} and the  coefficients $\tilde \delta_p$.

\smallskip

\begin{Proposition} \label{generalpropdeltaprime} 
Let $p>2$. Assume that we can find a   non-decreasing sequence of positive integers $(m_k)_{k \geq 1}$   such that $m_k =o ( 3^{2k/p} k^{-1}) $, as $k \rightarrow \infty$, in such a way that the following conditions are satisfied:
\begin{equation}\label{first-cond-delta-prime}
\sum_{k \geq k_0} 3^{k(p-2)/2} \left( \sum_{\ell \ge k}\delta_2'(m_\ell)m_{\ell+1}^{1/2}\right)^p <\infty \, , \quad  \qquad \sum_{k \geq k_0}  \left( \sum_{\ell \ge k}\delta_p'(m_\ell)m_{\ell+1}^{1-1/p}\right)^p<\infty  \, ,
\end{equation}
\begin{equation}\label{fourth-cond-delta-prime}
 \quad\sum_{\ell \ge k}\delta_2'(m_\ell)m_{\ell+1}^{1/2} 
%+\min_{m\ge 1} \big( m3^{k(2-p)/p}+\sum_{\ell \ge m}\delta_2'(\ell)/ 
%(\ell+1)^{1/2}\big) 
=o(3^{k(2-p)/2p}
/\sqrt{\log k})\, ,
\end{equation}
and there exists $r\in ]p,\infty[$, such that 
\begin{equation}\label{second-cond-delta-prime}
\sum_{k\ge 0}3^{k(p-r)/p}m_k^{(r-2)/2}\, <\infty\, ,
\end{equation}
and 
\begin{equation}\label{third-cond-delta-prime}
\sum_{j\ge 1} \frac{(\delta_p'(j))^{p/r}}{j^{1/r}}\, <\infty\, .
\end{equation}
Then, \eqref{series-variance} holds.
Moreover, if $\sigma>0$, one can redefine $(X_n)_{n \geq 1}$ without changing its distribution on a (richer) probability space on which 
there exist iid random variables $(N_i)_{i \geq 1}$ with common distribution ${\mathcal N} (0, \sigma^2)$, such that,
\beq \label{conclusiongeneralprop}
\left | \,   S_n - n \E (X_1)- \sum_{i=1}^n N_i \,  \right |= o\left (n^{1/p} \right ) \, \text{ ${\mathbb P}$-a.s.}
\eeq
\end{Proposition}

\subsection{Proof of Proposition \ref{generalpropdeltaprime}} 
\label{section-proof}

We consider a process $(X_n)_{n\ge 1}$ satisfying \eqref{functionofarandomiterates}, with stationary distribution 
$\pi$. 
We shall check that the assumptions of Proposition \ref{generalpropBLW} 
are satisfied.

Set $V_0:=(\varepsilon_0,W_0)$. It is not difficult to see that $\Vert   \E ( (X_{k}-\E(X_1))  | V_0 ) \Vert_2
\le \delta_2'(k)\le \delta_p'(k)$.  Now, since $r\ge p\ge 2$,
it follows from \eqref{third-cond-delta-prime}   that 
\beq\label{heyde-delta'}
\sum_{n\ge 1}\frac{\delta_2'(n)}{\sqrt n} <\infty\, .
\eeq
Then, the fact that \eqref{series-variance} holds 
follows from the fact that (see e.g. Lemma 22 of \cite{CDM})
\beq\label{ser-cov}
\sum_{k \geq 1}  |   {\rm Cov} (X_1, X_{k+1} )| \ll  \left (  \sum_{k \geq 0}  (k+1)^{-1/2} \Vert   \E ( (X_{k}-\E(X_1))  | V_0 ) \Vert_2  \right )^2< \infty \, .
\eeq

%Assume that $\sigma^2 >0$. 
We choose $M_k=3^{k/p}$. 
%Then, \eqref{condtruncature} follows from a standard computation. 
%Note that the sequence $(m_k)_{k \geq 0 }$ satisfies $m_k =o ( %3^{2k/p} k^{-1}) $, as $k \rightarrow \infty$.  
Since the $X_i$'s are in ${\mathbb L}^p$, it is easy to see that with this choice of $M_k$, condition  \eqref{condtruncature} is satisfied (it suffices to write that  $\E (|g_k (X_1)|)\leq   \E (|X_1|{\mathbf 1}_{|X_1| >M_k})$ and to use Fubini's Theorem).  

We shall check the condition \eqref{condmkdependence} with 
$\alpha=p$. To do so, we apply the Rosenthal-type inequality
given in 
Proposition \ref{prop-rosenthal} of the appendix, to the process 
$(\tilde X_{k,\ell+3^{k-1}}-X_{k,\ell+3^{k-1}})_{\ell \ge 1}$, with the choice 
$\eta_0 =  ( \varepsilon_{3^{k-1}},  \ldots, \varepsilon_{1+3^{k-1}-m_k}, W_{3^{k-1} -m_k}  )$
and 
for every $\ell\ge 1$, $\eta_\ell=\varepsilon_{\ell+3^{k-1}}$. 
We have to bound, for $q\ge 1$, the coefficients $\delta^*_{q,k}(n)$ used in Proposition 
\ref{prop-rosenthal}. When $1\le n\le m_k$, we use the bound 
$\delta^*_{q,k}(n)\le 2\|\tilde X_{k,3^{k-1}}-X_{k,3^{k-1}}\|_q
\le 2\delta'_q(m_k)$. When $n>m_k$ we notice that the contribution 
of  $(\tilde X_{k,\ell+3^{k-1}})_{\ell \ge 1}$ to $\delta^*_{q,k}(n)$ is null, so that 
$\delta^*_{q,k}(n)\le \delta_q'(n)$.

In particular we infer that,  for $k \geq k_0$,  

\begin{multline*}
\left \|\max_{1\le \ell \le 3^k-3^{k-1}}|\widetilde W_{k,\ell}-W_{k,\ell}|\,  \right \|_p \\
\ll \big(3^{k/2}\delta_2'(m_k)m_k^{1/2}+3^{k/p}m_k^{1-1/p}\delta_p'(m_k)
\big) 
+\left(3^{k/2}\sum_{j\ge m_k}\frac{\delta_2'(j)}{\sqrt j} +3^{k/p}\sum_{ 
j\ge m_k} \frac{\delta_p'(j)}{j^{1/p}}\right )\\
\ll  3^{k/2}\sum_{\ell \ge k}\delta_2'(m_\ell)m_{\ell+1}^{1/2}+
3^{k/p}\sum_{\ell \ge k} \delta_p'(m_\ell)m_{\ell+1}^{1-1/p}\, .
\end{multline*}
Hence,  \eqref{condmkdependence} holds with 
$\alpha=p$, since \eqref{first-cond-delta-prime} is satisfied.

%\[
%\sum_{k \geq k_0} 3^{k(p-2)/2} \Big( \sum_{\ell \ge k}\delta_2'(m_\ell)m_{\ell+1}^{1/2}\Big)^p <\infty \qquad {\rm and} \qquad \sum_{k \geq k_0}  \Big( \sum_{\ell \ge k}\delta_p'(m_\ell)m_{\ell+1}^{1-1/p}\Big)^p<\infty  \, ,
%]
%which is exactly \eqref{first-cond-delta-prime}.

We prove now that \eqref{condtoapplySakhanenko} holds for some $r>2$. 
We apply again Proposition \ref{prop-rosenthal}, but now to the process 
$(\tilde X_{k,\ell+3^{k-1}})_{1\le \ell \le 3m_k}$ and with the choice $\eta_0 =  ( \varepsilon_{3^{k-1}},  \ldots, \varepsilon_{1+3^{k-1}-m_k} )$
and 
for every $\ell\ge 1$, $\eta_\ell=\varepsilon_{\ell+3^{k-1}}$.

For every $q\ge 1$, denote by $\delta_{q,k}^*(n)$ the $n^{\rm th}$ coefficient $\delta^*$ associated with the above choice, and notice that 
$\delta_{q,k}^*(n)=0$ as soon as $n>m_k$. For every $q\ge 1$, denote by $\delta_{q,k}'(n)$ the $n^{\rm th}$ coefficient $\delta'$ 
associated with the process $(X_{k,\ell+3^{k-1}})_{\ell\ge 1}$. 
One can see that for every $n\ge 0$,  $\delta_{q,k}^*(n)\le \delta_{q,k}'(n)$.

For every $r\ge 2$, every $k\ge 1$, with $d_k$ the unique integer such that $2^{d_k-1}<3m_k\le 2^{d_k}$, Proposition \ref{prop-rosenthal} gives

\beq  \label{Wnormr-delta-prime}
 \left \Vert  \max_{1 \leq \ell \leq  3m_k }  
 \left |   
 {\widetilde W}_{k, \ell}    \right |\,  
 \right \Vert_r \ll 2^{d_k/2}\sum_{j=0}^{m_k} \delta'_{2,k}(j)/(j+1)^{1/2}
 +
2^{d_k/r}\sum_{j=0}^{m_k} \delta'_{r,k}(j)/(j+1)^{1/r}\,.
\eeq
Hence, \eqref{condtoapplySakhanenko} holds for some $r>2$, if 
\beq \label{inter-cond0}
\sum_{k\ge 0} 3^{k(p-r)/p}m_k^{(r-2)/2}<\infty \qquad \mbox{and} 
\qquad \sup_{k\ge 0} \sum_{j\ge 0} \delta'_{2,k}(j)/(j+1)^{1/2}< \infty\, ,
\eeq
and 
\begin{equation}\label{inter-cond}
\sum_{k\ge 0}3^{k(p-r)/p}\Big(\sum_{j\ge 0} \delta'_{r,k}(j)/(j+1)^{1/r}
\Big)^r<\infty\, .
\end{equation}
The first part of \eqref{inter-cond0}  is exactly \eqref{second-cond-delta-prime}. Moreover since  $\varphi_k$ is $1$-Lipschitz, we have $\delta_{2,k}'(n)\le  \delta_{2}'(n)\le \delta_p'(n)$.  Hence the second part  of  
\eqref{inter-cond0} holds for some $r>2$ as
soon as  
\eqref{third-cond-delta-prime} does.

It remains to prove \eqref{inter-cond}. By H\"older's inequality,
$$
\sum_{j\ge 0} \frac{\delta'_{r,k}(j)}{(j+1)^{1/r}} \le 
\left( \sum_{j\ge 0} \frac{(\delta_p'(j))^{p/r}}{(j+1)^{1/r}}\right)^{(r-1)/r}
\left( \sum_{j\ge 0}\frac{(\delta_{r,k}'(j))^r}{(j+1)} \frac{(j+1)^{(r-1)/r}}
{(\delta_p'(j))^{p(r-1)/r})}\right)^{1/r}\, .
$$
Taking into account \eqref{third-cond-delta-prime}, we see that 
\eqref{inter-cond}  holds as soon as 
$$
\sum_{j\ge 0}\frac{(\delta_p'(j))^{p(1-r)/r}}{(j+1)^{1/r}}
\sum_{k\ge 0} 3^{k(p-r)/p}(\delta_{r,k}'(j))^r<\infty\, .
$$
Using the fact that $\delta_{r,k}'(j) \leq 2 \Vert \varphi_k (X_0) \Vert_{r}$, that for every 
non negative random variable $Z$, $$\sum_{k\ge 0} 3^{k(p-r)/p}
\E((\varphi_k(Z))^r)\le C_{r,p}\E(Z^p)\, ,
$$ 
and \eqref{third-cond-delta-prime} again, 
we see that 
\eqref{inter-cond} holds.

To end the proof, it remains to prove that \eqref{cond2v_k} holds. 
Since $\sigma>0$, it follows from equation (65) of \cite{CDM} that 
 \eqref{cond2v_k} is satisfied as soon as
\beq \label{cond2v_kbis}
3^k (  \nu_k - \sigma^2)^2 = o(3^{2k/p} (\log k)^{-1})  \, , \, \mbox{ as $k \rightarrow \infty$}\, .  
\eeq
To prove \eqref{cond2v_kbis}, let us define, for $i \geq 0$, 
\[
{\tilde c}_{k,i}= {\rm cov} ({\tilde X}_{k, m_k+1} , {\tilde X}_{k,i+m_k+1} ) \, \text{ and } {\hat c}_{k,i} = {\rm cov} ({ X}_{k,0} , { X}_{k,i} ) \, .
\]
We have
\begin{equation}\label{est-var}
|\nu_k-\sigma^2|\le \left |\sum_{i=-m_k}^{m_k} \tilde c_{k,|i|}
-\sum_{i\in \Z}\hat c_{k,|i|}\right|+ \left|\sum_{i\in \Z} \hat c_{k,|i|}
-\sum_{i\in \Z} c_{|i|}\right|
\, .
\end{equation}
Arguing as in \cite{CDM} to obtain their equation (68), and making use of 
\eqref{heyde-delta'}, we see that 
$$
\left|\sum_{i=-m_k}^{m_k} \tilde c_{k,|i|}
-\sum_{i\in \Z}\hat c_{k,|i|}\right|\le C \Big ( \limsup_{j\to \infty}
j^{-1/2} \left \|\widetilde W_{k,j}-W_{k,j}
\right \|_2 +  \limsup_{j\to \infty}
j^{-1} \left \|\widetilde W_{k,j}-W_{k,j}
\right \|^2_2 \Big ) \, ,
$$
for some  $C>0$, independent of $k>0$.
Estimating the right-hand side    thanks to Proposition \ref{prop-rosenthal} 
with $p=2$, we infer that 
\begin{multline*}
 \left|\sum_{i=-m_k}^{m_k} \tilde c_{k,|i|}
-\sum_{i\in \Z}\hat c_{k,|i|}\right|\ll \left(1+\sum_{\ell\ge 0}\frac{\delta_2'(\ell)}{\sqrt {\ell+1}}\right)
\left(\sqrt{m_k}\delta_2'(m_k)+\sum_{\ell \ge m_k} 
\frac{\delta_2'(\ell)}{\sqrt {\ell+1}}\right)\\
\label{1st-term-var} \ll \left(\sum_{\ell\ge 0}\frac{\delta_2'(\ell)}{\sqrt {\ell+1}}\right )\sum_{\ell\ge k} \delta_2'(m_\ell) m_{\ell+1}^{1/2}\, . 
\end{multline*}
Hence by \eqref{fourth-cond-delta-prime} and \eqref{third-cond-delta-prime},
\begin{equation} \label{majcov1}
3^k \left|\sum_{i=-m_k}^{m_k} \tilde c_{k,|i|}
-\sum_{i\in \Z}\hat c_{k,|i|}\right|^2 = o(3^{2k/p}
(\log k)^{-1} )  \, . 
\end{equation}

Let now $c_i = {\rm Cov} (X_0,X_i) $ and note that (see Relation (3.54) in \cite{BLW14})
\[
\sup_{i  \geq 0} | {\hat c}_{k,i} -  c_i  | = o (3^{-k (p-2)/p}) \, .
\]
Let 
\[
\ell_k = 3^{k (p-2)/ (2p)} (\log k)^{-1/2} \, .
\]
It follows that
\[
 \left|\sum_{i\in \Z} \hat c_{k,|i|}
-\sum_{i\in \Z} c_{|i|}\right| \leq o ( \ell_k 3^{-k (p-2)/p}) + 2 \sum_{i > \ell_k} | c_{i} - { \hat c}_{k,i}| \, .
\]
Now
\[
| c_{i} - { \hat c}_{k,i}|  = |\cov ( X_0 - \varphi_k (X_0), X_i) + \cov (  \varphi_k (X_0), X_i - \varphi_k (X_i) )  | \, .
\]
Therefore
\begin{multline} \label{comp1sigma2avechat}
 \left|\sum_{i\in \Z} \hat c_{k,|i|}
-\sum_{i\in \Z} c_{|i|}\right| \leq o ( \ell_k 3^{-k (p-2)/p}) + 2 \sum_{i > \ell_k} |\cov ( X_0 - \varphi_k (X_0), X_i)  | \\ +  2 \sum_{i > \ell_k} |\cov (  \varphi_k (X_0), X_i - \varphi_k (X_i) )  |  \, .
\end{multline}
Let us first handle the series
\[
 \sum_{i > \ell_k} |\cov ( X_0 - \varphi_k (X_0), X_i)  | \, .
\]
Set $g_k (x) = x - \varphi_k(x)$. Applying Lemma 22 of \cite{CDM} and using the fact that $(W_k)_{k \geq 0}$ is a Markov chain, we infer that  
\begin{multline*}
 \sum_{i > \ell_k} |\cov ( X_0 - \varphi_k (X_0), X_i)   | \\
%  \ll  \sum_{j=0}^{\infty} \Vert P_{0} ( g_k (X_j)) \Vert_2 \sum_{i \geq [2^{-1} (\ell_k+ j)] +1 } 
%i^{-1/2} \Vert \E (X_i | V_0) \Vert_2 \\ 
 \ll \left (\sum_{i \geq [\ell_k/2] } 
i^{-1/2} \Vert \E (X_i | V_0) - \E (X_i)  \Vert_2\right ) \sum_{j=0}^{\infty} 
(j+1)^{-1/2}\Vert  \E (  g_k (X_j) | V_0) -  \E (  g_k (X_j))  \Vert_2\, ,
\end{multline*}
%where $P_0 (\cdot) = \E ( \cdot | {\mathcal F}_0 ) -  \E ( \cdot | {\mathcal F}_{-1} )$ with ${\mathcal F}_k= \sigma (X_\ell , \ell \leq k )$. 
We shall now use the following  estimate, to be proved at the end of this subsection: let $\varphi_M(x) = (x \wedge M)\vee (-M)$ and $g_M = x- \varphi_M(x)$, then
\beq \label{otherboundforcarrebis}
\Vert \E ( g_M (X_n) | V_0) - \E (g_M (X_n) )\Vert_2 \ll  \frac{1}{M^{(p-2)/2}   } \big (\delta_p'(n) \big )^{p/(2 (p-1))}   \, .
\eeq
Taking into account  \eqref{otherboundforcarrebis} and the fact that $\Vert \E (X_i | V_0) - \E (X_i)  \Vert_2 \leq \delta_2'(i) \leq \delta_p'(i)$, it follows that  
\begin{equation*}
 \sum_{i > \ell_k} |\cov ( X_0 - \varphi_k (X_0), X_i)   | \ll  
\frac1{M_k^{(p-2)/2}}\left( \sum_{i\ge [\ell_k/2]} i^{-1/2}\delta_p'(i)  
\right)\sum_{j=0}^{\infty} (j+1)^{-1/2}(\delta_p'(j))^{p/(2(p-1))}
\, .
\end{equation*}
Now, using \eqref{third-cond-delta-prime} and the fact that 
$(\delta_p'(j))_{j\ge 0}$ is non increasing, we see that 
$\delta_p'(j)=o(j^{(1-r)/p})$. In particular, $\sum_{j=0}^{\infty} (j+1)^{-1/2}(\delta_p'(j))^{p/(2(p-1))}<\infty$ and, since $r\ge p$,
$$
\sum_{i\ge [\ell_k/2]} i^{-1/2}\delta_p'(i) = O(\ell_k^{1/p-1/2})\, .
$$
Hence
\beq \label{majcov2}
 \sum_{i > \ell_k} |\cov ( X_0 - \varphi_k (X_0), X_i)   | \ll  \frac{\ell_k^{1/p-1/2} }{M_k^{(p-2)/2}}\, .
\eeq

Let us now handle  the series
\[
 \sum_{i > \ell_k} |\cov (  \varphi_k (X_0), X_i - \varphi_k (X_i) )     | \, .
\]
Applying again Lemma 22 of \cite{CDM} and taking into account the fact that $(W_k)_{k \geq 0}$ is a Markov chain,   we first infer that 
\begin{multline*}
 \sum_{i > \ell_k} |\cov (  \varphi_k (X_0), X_i - \varphi_k (X_i) )     | \\  
 \ll \sum_{\ell=0}^{\infty} 
(\ell+1)^{-1/2}\Vert  \E (  \varphi_k (X_\ell) | V_0) -  \E (  \varphi_k (X_\ell))  \Vert_2  \sum_{i \geq  [2^{-1} (\ell_k+\ell)] +1} 
i^{-1/2} \Vert \E (g_k ( X_i) | {\mathcal F}_0) - \E (g_k ( X_i) ) \Vert_2 \, .
%\\
% \ll \sum_{\ell=0}^{\infty} 
%(\ell+1)^{-1/2}\Vert  \E (  \varphi_k (X_\ell) | V_0) -  \E (  \varphi_k (X_\ell))  \Vert_2 \sum_{i \geq  [2^{-1} \ell_k] +1} 
%i^{-1/2} \Vert \E (g_k ( X_i) | V_0) - \E (g_k ( X_i) ) \Vert_2 \, .
\end{multline*}
Since $\varphi_k$ is $1$-Lipschitz, we have $\Vert  \E (  \varphi_k (X_\ell) | V_0) -  \E (  \varphi_k (X_\ell))  \Vert_2 \leq \delta_2'(\ell) \leq \delta_p'(\ell)$. Therefore, since by assumption, $\sum_{\ell=0}^{\infty} 
(\ell+1)^{-1/2} \delta_p'(\ell) < \infty$, 
\[
 \sum_{i > \ell_k} |\cov (  \varphi_k (X_0), X_i - \varphi_k (X_i) )     | \ll \sum_{i \geq  [2^{-1} \ell_k] +1} 
i^{-1/2} \Vert \E (g_k ( X_i) |V_0) - \E (g_k ( X_i) ) \Vert_2 \, .
\]
Using \eqref{otherboundforcarrebis}, the fact that $\delta_p'(i)=o(i^{(1-r)/p})$ and that $r>p$, it follows
\begin{multline} \label{series2}
 \sum_{i > \ell_k} |\cov (  \varphi_k (X_0), X_i - \varphi_k (X_i) )     |  \ll  \frac1{M_k^{(p-2)/2}} \sum_{i \geq  [2^{-1} \ell_k] +1} 
i^{-1/2} (\delta_p'(i))^{p/(2(p-1))} \\
 \ll 
\frac{\ell_k^{(p-r)/(2(p-1))}}{M_k^{(p-2)/2}} \, .
\end{multline}
Starting from \eqref{comp1sigma2avechat} and using \eqref{majcov2} and \eqref{series2}, we get 
\begin{equation} \label{majcov3}
 \left|\sum_{i\in \Z} \hat c_{k,|i|}
-\sum_{i\in \Z} c_{|i|}\right| \ll o ( \ell_k 3^{-k (p-2)/p}) + \frac{\ell_k^{1/p-1/2} }{M_k^{(p-2)/2}} +  \frac{\ell_k^{(p-r)/(2(p-1))}}{M_k^{(p-2)/2}}  \, .
\end{equation}
Starting from \eqref{est-var} and taking into account \eqref{majcov1} and \eqref{majcov3}, the condition \eqref{cond2v_kbis} is satisfied (since $M_k =3^{k/p}$, $\ell_k = 3^{k (p-2)/ (2p)} (\log k)^{-1/2}$ and $r>p>2$). 
\hfill $\square$

\medskip

\noindent {\bf Proof of \eqref{otherboundforcarrebis}.} We start by noticing that, for any $\varepsilon >0$, 
\begin{multline} \label{newestimate1}
\Vert \E ( g_M (X_n) | V_0) - \E (g_M (X_n) \Vert_2^2 = \int \big | \E (g_M (X_{n,x}))  - \int \E (g_M (X_{n,y})) \nu (dy)  \big |^2  \nu (dx) \\
\leq \iint \big | \E (g_M (X_{n,x}) - g_M (X_{n,y}) ) \big |^2 \nu (dx)  \nu (dy)  \\
\leq 2 \iint \big | \E \big (  (g_M (X_{n,x}) - g_M (X_{n,y}) ) {\bf 1}_{\{|g_M (X_{n,x}) - g_M (X_{n,y})| \leq \varepsilon \} } \big ) \big |^2  \nu (dx)  \nu (dy)  \\
+ 2  \iint \big | \E \big (  (g_M (X_{n,x}) - g_M (X_{n,y}) ) {\bf 1}_{\{|g_M (X_{n,x}) - g_M (X_{n,y})| > \varepsilon \} } \big ) \big |^2  \nu (dx)  \nu (dy) \, .
\end{multline}
Now, 
\begin{multline}  \label{newestimate2}
 \iint \big | \E \big (  (g_M (X_{n,x}) - g_M (X_{n,y}) ) {\bf 1}_{\{|g_M (X_{n,x}) - g_M (X_{n,y})| \leq \varepsilon \} } \big ) \big |^2  \nu (dx)  \nu (dy) \\
 \leq \varepsilon   \iint\E \big | g_M (X_{n,x}) - g_M (X_{n,y}) \big |  \nu (dx) d \nu (dy)  \\
  \leq 2  \varepsilon  \E (|g_M (X_n)|)   \leq 2  \varepsilon  \E (|X_n| {\bf 1}_{\{|X_n|  > M \} })   \leq 2  \varepsilon M^{1 -p}  \E (| X_1)|^p)   \, .
\end{multline}
On the other hand,
\begin{multline}  \label{newestimate3}
 \iint \big | \E \big (  (g_M (X_{n,x}) - g_M (X_{n,y}) ) {\bf 1}_{\{|g_M (X_{n,x}) - g_M (X_{n,y})| > \varepsilon \} } \big ) \big |^2  \nu (dx)  \nu (dy)  \\
 \leq  \iint \E \big (  (g_M (X_{n,x}) - g_M (X_{n,y}) )^2 {\bf 1}_{\{|g_M (X_{n,x}) - g_M (X_{n,y})| > \varepsilon \} } \big )   \nu (dx)  \nu (dy) \\
 \leq \varepsilon^{2-p}   \iint\E \big | g_M (X_{n,x}) - g_M (X_{n,y}) \big |^p  \nu (dx)  \nu (dy)  \\
 \leq 2^{p-1} \varepsilon^{2-p}   \iint\E \big | X_{n,x} - X_{n,y} \big |^p \nu (dx)  \nu (dy)  
 + 2^{p-1} \varepsilon^{2-p}   \iint\E \big | \varphi_M (X_{n,x}) - \varphi_M (X_{n,y}) \big |^p \nu (dx)  \nu (dy)  \\
  \leq 2^p  \varepsilon^{2-p} (  \delta_p' (n) )^p \, .
\end{multline}
Starting from \eqref{newestimate1},  taking into account \eqref{newestimate2}  and \eqref{newestimate3},  and selecting $\varepsilon = M  \big (\delta_p'(n) \big )^{p/ (p-1)} $, the upper bound \eqref{otherboundforcarrebis} follows.

\subsection{Proof of Theorem \ref{KMTavecdelta'p}}

Assume that $\delta_p'(n)=O(n^{-\gamma})$, for some $\gamma>0$. 

We shall first assume that $\sigma >0$ and   apply Proposition \ref{generalpropdeltaprime}. We shall take 
$m_k=O(3^{k\beta/p})$ for some $2\ge \beta >0$. Hence, we have to find 
$r>p$, $\gamma>0$ and $2\ge \beta>0$ such that \eqref{first-cond-delta-prime}, 
\eqref{fourth-cond-delta-prime}, \eqref{second-cond-delta-prime} and 
\eqref{third-cond-delta-prime} hold. One easily sees that the condition 
\eqref{first-cond-delta-prime} holds provided that 
$$
p-2< \beta(2\gamma-1) \qquad {\rm and}\qquad \gamma>1-1/p\, ;
$$
the condition \eqref{fourth-cond-delta-prime} holds provided that 
\beq\label{cond-1}
p-2< \beta(2\gamma-1)\, ;
\eeq
the condition \eqref{second-cond-delta-prime} holds provided that 
\beq\label{cond-2}
r<1+\gamma p\, ;
\eeq
and the condition  \eqref{third-cond-delta-prime} holds provided that 
\beq\label{cond-3}
\beta(r-2)<2(r-p)\, .
\eeq
Notice that the condition $p-2< \beta(2\gamma-1)$ appears twice, that 
\eqref{cond-3} implies that $\beta\le 2$, and that 
the condition $\gamma>1-1/p$ is realized 
as soon as 
$r<1+\gamma p$. 

Hence, we have to find $\gamma,\beta>0$ and $r>p$ such that 
\eqref{cond-1}, \eqref{cond-2} and \eqref{cond-3} hold.  In particular, one has to find $\beta>0$ such that 
$(p-2)/(2\gamma -1)<\beta < 2(r-2)/(r-p)$, which is possible as soon as 
$(p-2)/(2\gamma -1)< 2(r-2)/(r-p)$. The latter condition is equivalent to
(provided that $4\gamma>p$, a condition to be checked at the end):
$$
r-2>\frac{(2p-4)(2\gamma-1)}{4\gamma-p} \, .
$$ 
Now, one can find $r>p$ satisfying the latter condition and 
\eqref{cond-2} as soon as $$2+\frac{(2p-4)(2\gamma-1)}{4\gamma-p}
<1+\gamma p\, , $$
 which is equivalent to 
$$
4p\gamma^2 -(p^2+4p-4)\gamma +3p-4 >0\, .
$$
Then, one finds that $$\gamma>\frac{(p-2)\sqrt{p^2+12p+4}\,+p^2+4p-4}{8p}$$
solves the problem. 

Hence, Proposition \ref{generalpropdeltaprime} applies and the Theorem is proved in the case $\sigma>0$, provided that our condition on $\gamma$ 
imples that $4\gamma >p$, but this  may be  easily checked.

\medskip

Assume now on that $\sigma=0$. Proceeding as in the proof of Theorem 1 of 
\cite{CDM} (see page 17), we see that it suffices to prove that 
$$\sum_{n\ge 1}\frac{\delta_p'(n)}{n^{2/p^2}}<\infty\, .$$
Hence, it is enough to prove that 
$$
\frac2{p^2}+\frac{(p-2)\sqrt{p^2+12p+4}\,+p^2+4p-4}{8p}>1\, ,
$$
which in turn is equivalent to (recall that $p>2$)
\begin{equation}\label{inequality}
p\sqrt{p^2+12p+4}>-(p^2-2p-8)\, .
\end{equation}
The  right-hand side of \eqref{inequality} is non-positive for $p\ge 1+2\sqrt 2$. 
Taking the squares of \eqref{inequality}, one can see 
that \eqref{inequality} holds for $p\in ]2,1+2\sqrt 2]$, which ends the proof of the theorem. \hfill $\square$

\section{Appendix}

\subsection{A Rosenthal-type inequality under dependence} \label{SecRI}

We shall state and prove our inequality in a more general framework 
than needed. It is not difficult to prove that the coefficients $(\delta_p^*(n))_{n\ge 0}$ defined 
by \eqref{gen-delta-prime} below are 
precisely the ones introduced in \eqref{def-delta-prime}, taking 
$\eta_0=W_0$ and for every $k\ge 1$, $\eta_k=\varepsilon_k$.

Let $(\eta_k)_{k\ge 0}$ be independent random variables (not necessarily identically distributed) and define for every $k\ge 1$,
$\G_{k,-1/2}:=\{\emptyset, \Omega)$ and $\G_{k,0}=\sigma\{\eta_k\}$ and for every $k\ge 2$ and every 
$0\le \ell\le k-1$ , 
$\G_{k,\ell}:= \sigma\{\eta_k, \eta_{k-1}, \ldots, \eta_{k-\ell}\}$.

Let $(X_n)_{n\ge 1}$ be a process given by $X_n:=f_n(\eta_n,\eta_{n-1},\ldots, \eta_0)$, for $n\ge 1$, where $f_n$ 
is a real-valued measurable function. Assume that for every $n\ge 1$, 
$\E(|X_n|)<\infty$ and that $\E(X_n)=0$. We want to prove a Rosenthal-type inequality for $S_n:= X_1+\cdots +X_n$, $n\ge 1$.

We shall need the following measure of dependence. Let $(\eta'_k)_{k\ge 0}$ 
be an independent copy of $(\eta_k)_{k\ge 0}$. For every 
$k\ge m+1$ and every $m\ge 0$, set $X_{k,m}'=f_k(\eta_k,\ldots ,\eta_{k-m},\eta'_{k-m-1}, 
\ldots ,\eta'_0)$ and, then, for every $n\ge 1$,
\begin{equation}\label{gen-delta-prime}
\delta_p^*(n):= \sup_{m\ge n-1}\sup_{k\ge m} \|X_k-X'_{k,m}\|_p\, .
\end{equation}
Define also $\delta_p^*(0)=\sup_{k\ge 0}\|X_k\|_p$.

For every $\ell\ge 1$, set $T_{0,\ell}:=\E(X_\ell|\G_{\ell,0})$. 
 For every  $d\ge 0$, every $0\le k\le d$ and every $1\le \ell 
\le 2^{d-k}$, set 

\begin{equation}\label{def-U}
U_{k,\ell}:= \sum_{j=(\ell-1)2^k+1}^{\ell2^k} \big(X_j-\E(X_j|\G_{\ell 2^k ,2^k-1})\big) = \sum_{j=(\ell-1)2^k+1}^{\ell2^k} 
\big(X_j-\E(X_j|\G_{j ,j-(\ell-1) 2^k-1 })\big)\, .
\end{equation}
%=\sum_{j=(\ell-1)2^k+1}^{\ell2^k} X_j
%-\E\Big( \sum_{j=(\ell-1)2^k+1}^{\ell2^k} X_j\, \Big| \, 
%\G_{\ell 2^k ,2^k-1}\Big)
and 
\begin{equation*}
T_{k+1,\ell} =\E\left( \big(U_{k,2\ell-1}+U_{k,2\ell}\big)|\G_{\ell 2^{k+1} ,2^{k+1}-1}\right)
\, .
\end{equation*}

\medskip

\begin{Proposition}\label{prop-rosenthal}
For every $d\ge 0$, we have 
\begin{equation}\label{pointwise-rosenthal}
\max_{1\le n\le 2^d}|S_n|\le \sum_{k=0}^d \max_{1\le \ell\le 2^{d-k}}
\left|U_{k,\ell}\right| 
+\sum_{k=0}^d \max_{1\le m\le 2^{d-k}}
\left |\sum_{\ell=1}^{m} 
T_{k,\ell}\right|\, .
\end{equation}
In particular, if $X_n\in {\mathbb L}^p$, for every $n\ge 1$ and some $p\ge 2$, we have 
\begin{gather}\label{rosenthal}
\left \|\max_{1\le n\le 2^d}|S_n|\, \right \|_p\le \sum_{k=0}^{d} 
\left (\sum_{\ell=1}^{2^{d-k}} \|U_{k,\ell}\|_p^p\right )^{1/p}
+C_p \sum_{k=0}^d 
\left( \left (\sum_{\ell=1}^{2^{d-k}} \| T_{k,\ell}\|_2^2\right )^{1/2}+\left (\sum_{\ell=1}^{2^{d-k}} \|\ T_{k,\ell}\|_p^p\right )^{1/p}\right)\\
\label{rosenthal-delta}
\le C'_p 2^{d/2}\sum_{j=0}^{2^d} \frac{\delta_2^*(j)}{(j+1)^{1/2}}
\, +\, 
C''_p2^{d/p}\sum_{j=0}^{2^d} \frac{\delta_p^*(j)}{(j+1)^{1/p}}\,,
\end{gather}
where $C_p$  is the best constant in the Rosenthal inequality for 
independent random variables,
$C'_p=\frac{C_p 2^{3/2}}{\sqrt 2-1}$ and $C''_p=\frac{2^{1+1/p}(C_p+1)}{2^{1/p}-1}$. 
\end{Proposition}

\noindent {\bf Proof.} The proof is done  by induction on $d\ge 0$. 
The case where $d=0$ follows from the decomposition 
$$
X_n=\big( X_n-\E(X_n|\sigma\{\varepsilon_n\})\big) + \E(X_n|\sigma\{\varepsilon_n\})\, .
$$

Assume now that $\eqref{pointwise-rosenthal}$ holds for some $d\ge 0$. Let us prove that it holds for $d+1$. 

For every $n\ge 1$, we have $S_n=\sum_{k=1}^n \big(X_k-\E(X_k|\G_{k,0})\big)\, 
+\sum_{k=1}^n \E(X_k|\G_{k,0}):= R_n+\sum_{\ell=1}^n T_{0,\ell}$. Hence 
\begin{equation*}
\max_{1\le n\le 2^{d+1}}|S_n|\le \max_{1\le n\le 2^{d+1}} |R_n|+ \max_{1\le n\le 2^{d+1}}\left|\sum_{\ell=1}^n T_{0,\ell}\right|\, .
\end{equation*}

Using that for every $m\ge 1$, we have $|R_{2m+1}|\le 
|R_{2m}|+|X_{2m+1}-\E(X_{2m+1}|\G_{2m+1,0})|$, we infer that 
\begin{equation}
\label{simple-dec}\max_{1\le n\le 2^{d+1}}|S_n| \le \max_{1\le n\le 2^{d}} |R_{2n}|+ \max_{1\le k\le 2^{d+1}}|U_{0,k}|+\max_{1\le n\le 
2^{d+1}}\left|\sum_{\ell=1}^n T_{0,\ell}\right|\, .
\end{equation}
We shall use our induction hypothesis to handle the first term in the 
right-hand side of \eqref{simple-dec}.

For every $m\ge 1$, let 
$$\tilde X_m:= X_{2m-1}-\E(X_{2m-1}|\G_{2m-1,0})+ X_{2m}-
\E(X_{2m}|\G_{2m,0})\,  \text{ and } \, \tilde S_m:=R_{2m} \, .$$ 
Set $\tilde \eta_m:=(\eta_{2m},\eta_{2m-1})$ for  $m \geq 1$ and  $\tilde \eta_0:=\eta_0$. 
Let also $\tilde \G_{k,\ell}:= \sigma\{\tilde \eta_k,\ldots \tilde 
\eta_{k-\ell}\}=\G_{2k,2\ell+1}$. Then, for every $(\ell-1)2^k+1\le j
\le \ell2^k$, using that $\G_{2j-1,0}\subset 
\G_{\ell2^{k+1},2^{k+1}-1}=\tilde\G_{\ell 2^k ,2^k-1}$ and that $\G_{2j,0}\subset 
\G_{\ell2^{k+1},2^{k+1}-1}$, we have

\begin{equation*}
\tilde X_j-\E(
\tilde X_j|\tilde\G_{\ell 2^k ,2^k-1})= (X_{2j-1}+X_{2j})- \E (
X_{2j-1}+X_{2j}|\G_{\ell2^{k+1},2^{k+1}-1})\, .
\end{equation*}
Hence, for every $k\ge 0$,
\begin{multline*}
\tilde U_{k,\ell}:= \sum_{j=(\ell-1)2^k+1}^{\ell2^k} \big(\tilde X_j-\E(
\tilde X_j|\tilde\G_{\ell 2^k ,2^k-1})\big)=\sum_{j=(\ell-1)2^{k+1}+1}^{\ell2^{k+1}} 
\big(  X_{j}-\E(X_j |\G_{\ell2^{k+1},2^{k+1}-1})\big)=U_{k+1,\ell}\, ,
\end{multline*}
and 
\begin{equation*}
\tilde T_{k+1,\ell} :=\E\Big( \big(\tilde U_{k,2\ell-1}+\tilde U_{k,2\ell}\big)|\tilde \G_{\ell 2^{k+1} ,2^{k+1}-1}\Big)
=\E\Big( \big(U_{k+1,2\ell-1}+U_{k+1,2\ell}\big)|\G_{\ell 2^{k+2} ,2^{k+2}-1}\Big)=T_{k+2,\ell}
\, .
\end{equation*}
Notice that we also have 
$$
\tilde T_{0,\ell}=T_{1,\ell}.
$$ 
Applying the  induction hypothesis, we infer that
\begin{multline*}
\max_{1\le n\le 2^d}|\tilde S_n|\le \sum_{k=0}^d \max_{1\le \ell\le 2^{d-k}}
\big|\tilde U_{k,\ell}\big| 
+\sum_{k=0}^d \max_{1\le m\le 2^{d-k}}\left |\sum_{\ell=1}^{m} 
\tilde T_{k,\ell}\right |
\\\le \sum_{k=0}^d \max_{1\le \ell\le 2^{d-k}}
\big|U_{k+1,\ell}\big| 
+\sum_{k=0}^d \max_{1\le m\le 2^{d-k}}\left |\sum_{\ell=1}^{m} 
T_{k+1,\ell}\right | \\ =
\sum_{k=1}^{d+1} \max_{1\le \ell\le 2^{d+1-k}}
\big|U_{k,\ell}\big| 
+\sum_{k=1}^{d+1} \max_{1\le m\le 2^{d+1-k}}\left |\sum_{\ell=1}^{m} 
T_{k,\ell}\right|\, ,
\end{multline*}
which, combined with \eqref{simple-dec} yields \eqref{pointwise-rosenthal} 
with $d+1$ in place of $d$.

To prove \eqref{rosenthal}, we notice that on the one hand, 
$(\max_{1\le \ell\le 2^{d-k}}
\|U_{k,\ell}\| )^p\le \sum_{ \ell=1}^{ 2^{d-k}}
\|U_{k,\ell}\|^p$ and on the other hand, for every $0\le k \le d$,
 the variables $(T_{k,\ell})_{1\le \ell \le 2^{d-k}}$ are independent. 
 Then, it is a direct consequence from \eqref{pointwise-rosenthal} and the Rosenthal inequality for independent variables. 
 
 Since
%To prove \eqref{rosenthal-delta0}, we notice that by \eqref{def-U}, if $q\in \{2,p\}$,
$\|T_{0,\ell}\|_q\le \delta^*_q(0)$ and  $\|T_{k+1,\ell}\|_q  
\le 2\sum_{j=1}^{2^k}\delta_q^*(j)$, we infer 
from \eqref{rosenthal} that 
\begin{equation}
\label{rosenthal-delta0}
\left \|\max_{1\le n\le 2^d}|S_n|\, \right \|_p 
\le 2C_p\sum_{k=0}^d  2^{(d-k)/2} \sum_{j=0}^{2^{k-1}
} \delta_2^*(j) +2^{d/p}\delta_p^*(0)+
(2C_p+1) \sum_{k=0}^d  2^{(d-k)/p} \sum_{j=1}^{2^k} \delta_p^*(j) \, ,
\end{equation}
and \eqref{rosenthal-delta} easily follows.

 \hfill $\square$

\medskip

\noindent {\bf Acknowledgement.} The first author is very thankful to 
the laboratories MAP5 and LAMA for their invitations, that made possible the present collaboration.

\end{document}